\documentclass[12pt]{article}
\usepackage{amsfonts,amssymb}
\newcommand{\II}{{\mathbb I}}
\newcommand{\ZZ}{{\mathbb Z}}
\newcommand{\eps}{{\varepsilon}}
\newcommand{\elp}{{{\cal A}_{q,p}(\widehat{sl}(N)_{c})}}
\newcommand{\sfrac}[2]{{\textstyle{\frac{#1}{#2}}}}
\newcommand{\half}{{\frac{1}{2}}}
\newcommand{\car}[2]{\left[\begin{array}{c}{#1}\\{#2}\end{array}\right]}
\newcommand{\finproof}{{\hfill \rule{5pt}{5pt}}}

\newtheorem{coro}{Corollary}
\newtheorem{lemm}{Lemma}
\newtheorem{thm}{Theorem}

\begin{document}
\newpage
\pagestyle{empty}
\setcounter{page}{0}
\vfill
\begin{center}

{\Large \textbf{Universal construction of ${\cal W}_{q,p}$ algebras}}

\vspace{10mm}

{\large J. Avan}

\vspace{4mm}

{\em LPTHE, CNRS-URA 280, Universit{\'e}s Paris VI/VII, France}

\vspace{7mm}

{\large L. Frappat, M. Rossi, P. Sorba}

\vspace{4mm}

{\em Laboratoire d'Annecy-le-Vieux de Physique Th{\'e}orique LAPTH, 
CNRS-URA 1436}

{\em LAPP, BP 110, F-74941 Annecy-le-Vieux Cedex, France}

\end{center}

\vfill
\vfill

\begin{abstract}
We present a direct construction of abstract generators for $q$-deformed 
${\cal W}_N$ algebras.  This procedure hinges upon a twisted trace 
formula for the elliptic algebra $\elp$ generalizing the previously 
known formulae for quantum groups.
\end{abstract}

\vfill
\vfill

\rightline{LAPTH-690/98}
\rightline{PAR-LPTHE 98-38}
\rightline{math.QA/9807048}
\rightline{July 1998}

\newpage
\pagestyle{plain}

%%%%%%%%%%%%%%%%%%%%%%%%%%%%%%%%%%%%%%%%%%%%%%%%%%%%%%%%%%%%%%%%%%%%%%%
\section{Introduction}

The connection between $q$-deformed Virasoro, and more generally ${\cal 
W}$ algebras, and elliptic quantum $\elp$ algebras, was investigated in 
our recent papers [1--4].  It was shown that $q$-deformed Virasoro and 
${\cal W}$ structures [5--7] were present inside the $\elp$ elliptic 
algebra \cite{FIJKMY,JKOS}; at the quantum level, when a particular 
relation existed between the central charge $c$, the elliptic nome $p$ 
and the deformation parameter $q$: $(-p^{\frac{1}{2}})^{NM} = q^{-c-N}$ 
for some integer $M$, and at the classical limit, obtained when setting 
an additional relation $p = q^{Nh}$ for some integer $h$.  In this way, 
one obtained directly a set of quantizations of the classical $q$-deformed 
Poisson algebras, interestingly distinct from the original quantization 
\cite{FF} obtained from explicit $q$-boson realizations.

The construction was achieved at the abstract level in that only the 
abstract algebraic relations for $\elp$, defined by the eight vertex 
model $R$-matrix \cite{Bax}, were used to derive the $q$-deformed 
structures.  It was assumed throughout the derivations that the initial 
formal series relations \cite{FIJKMY} were in fact extended to the level 
of analytic relations, thereby leading from \emph{one} single exchange 
relation for this generating operator functional of the algebras to an 
\emph{infinite $\ZZ$-labeled set} of exchange relations for the modes, 
depending upon the choice of a relevant series expansion in a 
crown-shaped sector for the ratio of spectral parameters in the elliptic 
structure function.

In our original approach \cite{AFRS3}, the extension of the construction 
to $sl(N)$ was achieved by defining the abstract generators of higher 
spin simply as shifted ordered products of the spin one generators $t(z) 
= {\mbox{Tr}} \left[ L^+(zq^{c/2}) (L^-(z))^{-1} \right]$.  In this 
respect, the first derivation cannot be considered as the $q$-deformed 
version of the ${\cal W}_N$ algebra construction \cite{Bais} which takes 
as generators combinations of the current algebra generators from which 
one then extracts $sl(N)$ scalar objects; the detailed study developed 
in \cite{Bais} allows to appreciate the successes and the difficulties 
of this approach.  The question of an universal construction of ${\cal 
W}_{q,p}$ algebras from elliptic algebras thus remained open, although 
the construction \cite{AFRS3} gave rise to perfectly consistent 
non-trivial algebraic structures, due to the shift in the spectral 
parameters which did no allow to interprete such ${\cal W}_{q,p}$ 
algebras as simply enveloping algebras of $Vir_{q}(sl(2))$.  In 
particular, the classical limit of our algebra ${\cal W}_{q,p}$ did lead 
to the original \cite{FR} classical $q$-${\cal W}$ Poisson algebra, 
characterizing the quantum structure as a genuine $q$-deformed ${\cal 
W}_N$ algebra.

We address here this question.  We shall rely on the algebraic structures 
derived from the properties of the operator ${\mathfrak T}(z) \equiv 
L^+(q^{\frac{c}{2}}z) \, (L^-(z))^{-1}$ which was the fundamental  
object in our previous derivation.

In a first part, we prove that ${\mathfrak T}(z)$ obeys an exchange 
relation of the type $R\,{\mathfrak T}\,R'\,{\mathfrak T} = 
R'\,{\mathfrak T}\,R\,{\mathfrak T}$.  Originally derived and discussed 
in \cite{RSTS,MF}, these exchange relations then lead us to define new 
surfaces in the $(p,q,c)$ space on which quantum, then classical, 
$q$-Virasoro algebras of the same type as in \cite{AFRS1} arise. 

The classical structures are the same as in \cite{AFRS1}.  The quantum 
structures by contrast are more general, for $N \geq 3$, than the 
original algebras derived in \cite{AFRS3}, which one recovers as 
particular cases.  One cannot however directly derive higher order 
generators from such an exchange algebra, contrary to the simpler case 
$RLL = LLR$ where a famous twisted trace formula exists \cite{Mail,ABB} 
to generate quantum commuting Hamiltonians.  But since the definition of 
the basic elliptic algebra involves two distinct $R$-matrices as $RLL = 
LLR^*$, one cannot apply \cite{Mail,ABB} to it either.

In a second part, we show how to define a suitable twisted trace 
formula, involving ${\cal L}(z) \equiv L^+(q^{\frac{c}{2}}z)^t \, 
(L^-(z)^{-1})^t$ and the $R$-matrix of the algebra $\elp$, and leading 
to closed exchange algebras of the quantum ${\cal W}_{q,p}$ type, for 
the generalized relation $(-p^{\frac{1}{2}})^{n} = q^{-c-N}$ where $n$
is any integer, not necessarily multiple of $N$; then to 
classical $q$-${\cal W}$ Poisson algebras when $p = q^{Nh}$.

The quantum and classical algebras ${\cal W}_{q,p}[sl(N)]$ thus 
constructed contain in particular the structures obtained in 
\cite{AFRS3} by using the shifted products.  As emphasized, we now have 
a really \emph{universal} construction of the higher $q$-deformed 
currents for these ${\cal W}_{q,p}$ algebras.  In addition, we have 
obtained the elliptic algebra version of the twisted trace construction, 
bypassing thereby the difficulty generated a priori by the asymmetrical 
Yang--Baxter relation $RLL = LLR^*$ ($R^* \ne R$) defining the elliptic 
algebra.

%%%%%%%%%%%%%%%%%%%%%%%%%%%%%%%%%%%%%%%%%%%%%%%%%%%%%%%%%%%%%%%%%%%%%%%
\section{The elliptic algebra $\elp$ and the algebras  $Vir_{q,p}(sl(N))$}
\setcounter{equation}{0}

We first recall the most important notations, definitions and results 
concerning the elliptic algebras and their subsequent $q$-deformed 
Virasoro and ${\cal W}_N$ algebras.

\subsection{Definition of the elliptic $R$-matrix}

We start by defining the $R$-matrix of the $\ZZ_{N}$-vertex model 
($\ZZ_N$ is the congruence ring modulo $N$) \cite{Bax,Bela}:
\begin{equation}
\widetilde R(z,q,p) = z^{2/N-2} \frac{1}{\kappa(z^2)} 
\frac{\vartheta\car{\sfrac{1}{2}}{\sfrac{1}{2}}(\zeta,\tau)} 
{\vartheta\car{\sfrac{1}{2}}{\sfrac{1}{2}}(\xi+\zeta,\tau)} \,\, 
\sum_{(\alpha_1,\alpha_2)\in\ZZ_N\times\ZZ_N} 
W_{(\alpha_1,\alpha_2)}(\xi,\zeta,\tau) \,\, I_{(\alpha_1,\alpha_2)} 
\otimes I_{(\alpha_1,\alpha_2)}^{-1} \,,
\label{eq21}
\end{equation}
where the variables $z,q,p$ are related to the $\xi,\zeta,\tau$ 
variables by
\begin{equation}
z=e^{i\pi\xi} \,,\qquad q=e^{i\pi\zeta} \,,\qquad p=e^{2i\pi\tau} \,.
\end{equation}
$\vartheta$ are the Jacobi theta functions with rational 
characteristics $(\gamma_1,\gamma_2) \in \ZZ/N \times \ZZ/N$:
\begin{equation}
\vartheta\car{\gamma_1}{\gamma_2}(\xi,\tau) = \sum_{m \in \ZZ}
\exp\Big(i\pi(m+\gamma_1)^2\tau + 2i\pi(m+\gamma_1)(\xi+\gamma_2)
\Big) \,.
\end{equation}
The normalization factor is given by:
\begin{equation}
\frac{1}{\kappa(z^2)} = \frac{(q^{2N}z^{-2};p,q^{2N})_\infty
\, (q^2z^2;p,q^{2N})_\infty \, (pz^{-2};p,q^{2N})_\infty \,
(pq^{2N-2}z^2;p,q^{2N})_\infty} {(q^{2N}z^2;p,q^{2N})_\infty
\, (q^2z^{-2};p,q^{2N})_\infty \, (pz^2;p,q^{2N})_\infty \,
(pq^{2N-2}z^{-2};p,q^{2N})_\infty} \,.
\label{eq24}
\end{equation}
The functions $W_{(\alpha_1,\alpha_2)}$ are defined as follows:
\begin{equation}
W_{(\alpha_1,\alpha_2)}(\xi,\zeta,\tau) = \frac{1}{N} \,\, 
\frac{\vartheta\car{\sfrac{1}{2}+\alpha_1/N}
{\sfrac{1}{2}+\alpha_2/N}(\xi+\zeta/N,\tau)} 
{\vartheta\car{\sfrac{1}{2}+\alpha_1/N}
{\sfrac{1}{2}+\alpha_2/N}(\zeta/N,\tau)} \,,
\end{equation}
and the matrices $I_{(\alpha_1,\alpha_2)}$ by:
\begin{equation}
I_{(\alpha_1,\alpha_2)} = g^{\alpha_{2}} \, h^{\alpha_{1}} \,,
\end{equation}
where $g_{ij} = \omega^{i}\delta_{ij}$, $h_{ij} = \delta_{i+1,j}$ are $N 
\times N$ matrices (the addition of indices being understood modulo $N$) 
and $\omega = e^{2i\pi/N}$.  \\
The $R$-matrix $\widetilde R$ is $\ZZ_N$-symmetric:
\begin{equation}
\widetilde R_{c+s\,,\,d+s}^{a+s\,,\,b+s} = \widetilde 
R_{c\,,\,d}^{a\,,\,b} \qquad a,b,c,d,s \in \ZZ_N \,.
\end{equation}
%$p$ and $q$ restricted to $\vert q \vert < 1$, $\vert p \vert < 1$.
We introduce the ``gauge-transformed'' matrix:
\begin{equation}
R(z,q,p) = (g^{\frac 12} \otimes g^{\frac 12}) \widetilde
R(z,q,p) (g^{-\frac 12} \otimes g^{-\frac 12})
\end{equation}
It satisfies the following properties: \\
-- Yang--Baxter equation:
\begin{equation}
R_{12}(z) \, R_{13}(w) \, R_{23}(w/z) =
R_{23}(w/z) \, R_{13}(w) \, R_{12}(z) \,,
\label{eq29}
\end{equation}
-- unitarity: 
\begin{equation}
R_{12}(z) R_{21}(z^{-1}) = 1 \,,
\end{equation}
-- crossing symmetry: 
\begin{equation}
R_{12}(z)^{t_2} R_{21}(q^{-N}z^{-1})^{t_2} = 1 \,,
\end{equation}
-- antisymmetry: 
\begin{equation}
R_{12}(-z) = \omega \, (g^{-1} \otimes \II) \, R_{12}(z)
\, (g \otimes \II) \,,
\end{equation}
-- quasi-periodicity: 
\begin{equation}
\widehat{R}_{12}(-p^{\frac{1}{2}}z) = 
(g^{\frac{1}{2}}hg^{\frac{1}{2}} \otimes \II)^{-1} \left 
(\widehat{R}_{21}(z^{-1})\right )^{-1} (g^{\frac{1}{2}}hg^{\frac{1}{2}} 
\otimes \II) \,, 
\label{eq213}
\end{equation}
where
\begin{equation}
\widehat{R}_{12}(x) = \tau_{N}(q^{1/2}x^{-1}) \, R_{12}(x) \,,
\label{eq214}
\end{equation}
and
\begin{equation}
\tau_N(z) = z^{\frac{2}{N}-2} \,
\frac{\Theta_{q^{2N}}(qz^2)}{\Theta_{q^{2N}}(qz^{-2})} \,.
\label{eq215}
\end{equation}
The function $\tau_N(z)$ is periodic with period $q^N$:
$\tau_N(q^Nz) = \tau_N(z)$ and satisfies $\tau_N(z^{-1}) =
\tau_N(z)^{-1}$.

\subsection{Definition of $\elp$}

We now recall the definition of the elliptic quantum algebra $\elp$
\cite{FIJKMY,JKOS}. It is an algebra of operators $L_{ij}(z) \equiv 
\sum_{n \in \ZZ} L_{ij}(n) \, z^{n}$ where $i,j \in \ZZ_{N}$:
\begin{equation}
L(z) = \left(\begin{array}{ccc} L_{11}(z) & \cdots &
L_{1N}(z) \cr \vdots && \vdots \cr L_{N1}(z) & \cdots &
L_{NN}(z) \cr \end{array}\right) \,.
\end{equation}
The $q$-determinant is given by ($\eps(\sigma)$ being the signature of 
the permutation $\sigma$):
\begin{equation}
q\mbox{-}\det L(z) \equiv \displaystyle \sum_{\sigma\in{\mathfrak S}_N}
\eps(\sigma) \prod_{i=1}^N L_{i,\sigma(i)}(z q^{i-N-1}) \,.
\end{equation}
${\cal A}_{q,p}(\widehat{gl}(N)_c)$ is defined by imposing the following 
constraints on the $L(z)$ generators:
\begin{equation}
\widehat{R}_{12}(z/w) \, L_1(z) \, L_2(w) = L_2(w) \, L_1(z) \,
\widehat{R}_{12}^{*}(z/w) \,,
\label{eq218}
\end{equation}
where $L_1(z) \equiv L(z) \otimes \II$, $L_2(z) \equiv \II \otimes L(z)$ 
and $\widehat R^{*}_{12}(z,q,p) \equiv \widehat 
R_{12}(z,q,p^*=pq^{-2c})$.  \\
The $q$-determinant is in the center of ${\cal 
A}_{q,p}(\widehat{gl}(N)_c)$ and one sets
\begin{equation}
{\cal A}_{q,p}(\widehat{sl}(N)_c) = {\cal A}_{q,p}(\widehat{gl}(N)_c)/ 
\langle q\mbox{-}\det L - q^{\frac{c}{2}} \rangle \,.
\end{equation}

\subsection{$Vir_{q,p}(sl(N))$ algebras from elliptic algebras}

It is here useful to introduce the following two matrices:
\begin{equation}
L^+(z) \equiv L(q^{\frac{c}{2}}z) \,, \qquad \qquad
L^-(z) \equiv (g^{\frac 12} h g^{\frac 12}) \,
L(-p^{\frac 12}z) \, (g^{\frac 12} h g^{\frac 12})^{-1} \,,
\label{eq220}
\end{equation}
which obey coupled exchange relations following from (\ref{eq218}) and
periodicity/unitarity properties of the matrices
$\widehat R_{12}$ and $\widehat R^{*}_{12}$:
\begin{eqnarray}
&& \widehat R_{12}(z/w) \, L^\pm_1(z) \,
L^\pm_2(w) = L^\pm_2(w) \, L^\pm_1(z) \,
\widehat R^{*}_{12}(z/w) \,, \nonumber \\
&& \widehat R_{12}(q^{\frac{c}{2}}z/w) \,
L^+_1(z) \, L^-_2(w) = L^-_2(w) \, L^+_1(z)
\, \widehat R^{*}_{12}(q^{-\frac{c}{2}}z/w) \,.
\label{eq221}
\end{eqnarray}

We now recall some of the main results of refs. [1--3]:
\begin{thm}\label{thmvir}
In the three-dimensional parameter space generated by $p,q,c$, one 
defines a two-dimen\-sional surface $\Sigma_{N,NM}$ for any integer $M 
\in \ZZ$ by the set of triplets $(p,q,c)$ connected by the relation 
$(-p^{\frac{1}{2}})^{NM} = q^{-c-N}$.  One defines the following 
operators in $\elp$:
\begin{equation}
t(z)  \equiv {\rm Tr}\Big( 
L^+(q^{\frac{c}{2}}z) \, L^-(z)^{-1} \Big) = {\rm Tr}\Big( 
L^+(q^{\frac{c}{2}}z)^t \, \widetilde L^-(z) \Big) \equiv {\rm Tr}({\cal 
L}(z)) \,,
\end{equation}
where $\widetilde L^-(z) \equiv (L^-(z)^{-1})^t$.  \\
1) On the surface $\Sigma_{N,NM}$, the operators $t(z)$ realize an 
exchange algebra with the generators $L(w)$ of $\elp$:
\begin{equation}
t(z) \, L(w) = F_{N}\Big(NM,\frac{w}{z}\Big) \, L(w) \, t(z)
\label{eq31}
\end{equation}
where
\begin{equation}
F_{N}(r,x) = \left\{ \begin{array}{ll} 
\displaystyle q^{2r(1-\frac{1}{N})} \prod_{k=0}^{r-1}
\frac{\Theta_{q^{2N}}(x^{-2}p^{-k}) \, \Theta_{q^{2N}}(x^2p^{k})}
{\Theta_{q^{2N}}(x^{-2}q^2p^{-k}) \, \Theta_{q^{2N}}(x^2q^2p^{k})}
& \mbox{for $r>0$} \,, \\ \\
\displaystyle q^{-2|r|(1-\frac{1}{N})} \prod_{k=1}^{|r|}
\frac{\Theta_{q^{2N}}(x^{-2}q^2p^{k}) \, \Theta_{q^{2N}}(x^2q^2p^{-k})}
{\Theta_{q^{2N}}(x^{-2}p^{k}) \, \Theta_{q^{2N}}(x^2p^{-k})}
& \mbox{for $r<0$} \,. \\
\end{array} \right.
\label{eq32}
\end{equation}
2) On the surface $\Sigma_{N,NM}$, $t(z)$ closes a quadratic subalgebra:
\begin{equation}
t(z)t(w) = {\cal Y}_{N}\Big(NM,\frac{w}{z}\Big) \, t(w)t(z)
\label{eq33}
\end{equation}
where
\begin{equation}
{\cal Y}_{N}(r,x) = \left\{ \begin{array}{ll}
\displaystyle \prod_{k=1}^{r} \frac{\Theta_{q^{2N}}^2(x^2 p^{-k}) \,
\Theta_{q^{2N}}(x^2 q^2 p^k) \, \Theta_{q^{2N}}(x^2 q^{-2} p^k)}
{\Theta_{q^{2N}}^2(x^2 p^k) \, \Theta_{q^{2N}}(x^2 q^2 p^{-k})
\, \Theta_{q^{2N}}(x^2 q^{-2} p^{-k})}
& \mbox{for $r>0$} \,, \\ \\
\displaystyle \prod_{k=0}^{|r|-1}
\frac{\Theta_{q^{2N}}^2(x^2 p^{-k}) \,
\Theta_{q^{2N}}(x^2 q^2 p^k) \, \Theta_{q^{2N}}(x^2 q^{-2} p^k)}
{\Theta_{q^{2N}}^2(x^2 p^k) \, \Theta_{q^{2N}}(x^2 q^2 p^{-k})
\, \Theta_{q^{2N}}(x^2 q^{-2} p^{-k})}
& \mbox{for $r<0$} \,. \\
\end{array} \right.
\label{eq34}
\end{equation}
3) In particular, at the critical level $c=-N$, the operators $t(z)$ 
lie in the center of $\elp$ and commute with each other.
\end{thm}

\section{The new $Vir_{q,p}(sl(N))$ algebras}
\setcounter{equation}{0}

We first prove the main basic result of this section.

\subsection{A generalized quadratic exchange algebra}

\begin{thm}\label{thmcurl}
The operators ${\mathfrak T}(z)$ defined by
\begin{equation}
{\mathfrak T}(z) \equiv L^+(q^{\frac{c}{2}}z) \, L^-(z)^{-1}
\label{eq222}
\end{equation}
satisfy the following exchange relation:
\begin{equation}
\widehat R_{12}(z/w) \, {\mathfrak T}_{1}(z) \widehat R_{21}(q^{c}w/z) 
\, {\mathfrak T}_{2}(w) = {\mathfrak T}_{2}(w) \, \widehat 
R_{12}(q^{c}z/w) \, {\mathfrak T}_{1}(z) \, \widehat R_{21}(w/z)
\label{eq223}
\end{equation}
\end{thm}

\medskip

\noindent 
\textbf{Proof:} One can derive from eqs.  (\ref{eq221}) further exchange 
relations between the operators $L^+$ and $\widetilde L^-$ \cite{AFRS3}.  
One has therefore:
\begin{eqnarray}
{\mathfrak T}_{1}(z) \, \widehat R_{21}(q^{c}w/z) \, {\mathfrak 
T}_{2}(w) &=& L^+_{1}(q^{\frac{c}{2}}z) \, \widetilde L^-_{1}(z)^{t_{1}} 
\, \widehat R_{21}(q^{c}w/z) \, L^+_{2}(q^{\frac{c}{2}}w) \, \widetilde 
L^-_{2}(w)^{t_{2}} \nonumber \\
&=& L^+_{1}(q^{\frac{c}{2}}z) \, L^+_{2}(q^{\frac{c}{2}}w) \, 
\widehat R_{21}^{*}(w/z) \, \widetilde L^-_{1}(z)^{t_{1}} \, 
\widetilde L^-_{2}(w)^{t_{2}} \nonumber \\
&=& \widehat R_{12}^{-1}(z/w) \, L^+_{2}(q^{\frac{c}{2}}w) \, 
L^+_{1}(q^{\frac{c}{2}}z) \, \widehat R_{12}^{*}(z/w) \,
\widehat R_{21}^{*}(w/z) \, \widetilde L^-_{1}(z)^{t_{1}} \, 
\widetilde L^-_{2}(w)^{t_{2}} \nonumber \\
&=& T \, \widehat R_{12}^{-1}(z/w) \, L^+_{2}(q^{\frac{c}{2}}w) \, 
L^+_{1}(q^{\frac{c}{2}}z) \, \widetilde L^-_{1}(z)^{t_{1}} \, 
\widetilde L^-_{2}(w)^{t_{2}} \nonumber \\
&=& T \, \widehat R_{12}^{-1}(z/w) \, L^+_{2}(q^{\frac{c}{2}}w) \, 
L^+_{1}(q^{\frac{c}{2}}z) \, \widehat R_{12}^{*}(z/w) \, 
\widetilde L^-_{2}(w)^{t_{2}} \, \widetilde L^-_{1}(z)^{t_{1}} \, 
\widehat R_{12}^{-1}(z/w) \nonumber \\
&=& \widehat R_{12}^{-1}(z/w) \, L^+_{2}(q^{\frac{c}{2}}w) \, 
\widetilde L^-_{2}(w)^{t_{2}} \, \widehat R_{12}(q^{c}z/w) \, 
L^+_{1}(q^{\frac{c}{2}}z) \, \widetilde L^-_{1}(z)^{t_{1}} \, 
\widehat R_{21}(w/z) \nonumber \\
&& \label{eq224}
\end{eqnarray}
where $T$ stands for $\tau_{N}(q^{\half}w/z)\tau_{N}(q^{\half}z/w)$ 
and we used the relations $\widehat R_{12}^{*}(z/w) \, \widehat 
R_{21}^{*}(w/z) = T$ and $T \, \widehat R_{12}^{-1}(z/w) = \widehat 
R_{21}(w/z)$.  \\
Multiplying the last equality of (\ref{eq224}) by $\widehat R_{12}(z/w)$ 
on the left, one obtains the desired equation (\ref{eq223}). 
 
\finproof

\subsection{An alternative construction: new surfaces in $\elp$}

Mixed exchange relations of the type described in Theorem 2 were 
considered in \cite{RSTS,MF}.  It was then shown in \cite{MF} that, 
provided that a $c$-number matrix exist such that $\widehat R_{12}(z/w) 
\, \gamma_1 \widehat R_{21}(q^{c}w/z) \, \gamma_2 = \gamma_2 \, \widehat 
R_{12}(q^{c}z/w) \, \gamma_1 \, \widehat R_{21}(w/z)$, one may construct 
commuting generators defined as $Q \equiv \mbox{tr}(\tilde \gamma^t 
{\mathfrak T})$ where $\tilde \gamma$ is matrix dual to $\gamma$.  This, 
together with the properties of quasi-periodicity and unitarity of the 
$R$-matrix, leads us to consider the following operator:
\begin{equation}
t(z) \equiv {\mbox{Tr}} \Big[ a^{-n} \, {\mathfrak T}(z) \Big] = 
{\mbox{Tr}} \Big[ a^{-n} L^+(q^{c/2}z) \widetilde L^-(z)^t \Big] \,,
\label{eq91}
\end{equation}
where $n\in \ZZ$ and the matrix $a$ is given by:
\begin{equation}
a = g^{1/2} h g^{1/2} \,.
\label{eq92}
\end{equation}

By analogy with the construction in Theorem 1, we expect that the 
mechanism of construction of directly commuting Hamiltonians in 
\cite{MF} will here turn into a two-step procedure, with a first 
constraint on $p, q, c$ leading to a closed exchange algebra and a 
second constraint leading to commuting operators and a subsequent 
Poisson structure.

Indeed we first establish the exchange properties of (\ref{eq91}) on the 
surfaces $\Sigma_{N,n}$ of the three-dimensional space of parameters 
$p$, $q$, $c$, given by the equation ($n \in \ZZ$, $n \ne 0$):
\begin{equation}
(-p^{1/2})^n = q^{-c-N} \,.
\label{eq93}
\end{equation}

\medskip

It is relevant for later purposes to rewrite the operator $t(z)$ as 
$t(z) = {\mbox{Tr}} \Big[ {\cal L}^{(n)}(z) \Big]$ where ${\cal 
L}^{(n)}(z)$ is defined by
\begin{equation}
{\cal L}^{(n)}(z) = \left( a^{-n} L^+(q^{c/2}z) \right)^{t} \widetilde 
L^-(z) \,.
\label{eq925}
\end{equation}

We prove the following lemma.
\begin{lemm}
On the surfaces $\Sigma_{N,n}$, the operators ${\cal L}^{(n)}(z)$ 
defined by (\ref{eq925}) have the following exchange properties with the 
generators $L(w)$ of $\elp$:
\begin{equation}
{\cal L}_1^{(n)}(z) \, L_2(w)= F_{N}\left(n,\frac{w}{z}\right) L_2(w) \left( 
\widehat R^{*}_{21}(q^{-c}w/z)^{-1} \right)^{t_1} {\cal L}_1^{(n)}(z) \, 
\widehat R^{*}_{21}(q^{-c}w/z)^{t_1} \,.  \label{eq926}
\end{equation}
\end{lemm}

\medskip

\noindent \textbf{Proof:} It is easier to formulate the proof in terms 
of $L^+(w)$.  One has:
\begin{equation}
{\cal L}_1^{(n)}(z) \, L_2^+(w) = L_1^+(zq^{c/2})^{t_1} \, 
(a_1^{-n})^{t_1} \, \widetilde L_1^-(z) \, L_2^+(w) \,.
\label{eq95}
\end{equation}
To exchange $t(z)$ with $L^+(w)$, we need the following exchange 
relations, coming directly from (\ref{eq221}):
\begin{eqnarray}
\left(\widehat R_{21}(q^{\frac{c}{2}}w/z)^{t_1}\right)^{-1} \, L^+_2(w) 
\, \widetilde L^-_1(z) &=& \widetilde L^-_1(z) \, L^+_2(w) \, 
\left({\widehat R_{21}^*(q^{-\frac{c}{2}}w/z)}^{t_1}\right)^{-1} \,, 
\label{eq96} \\
L_1^+(z)^{t_1} \, \widehat R_{12}(z/w)^{t_1} \, L_2^+(w) &=& L_2^+(w) \, 
\widehat R_{12}^*(z/w)^{t_1} \, L_1^+(z)^{t_1} \label{eq97} \, .
\end{eqnarray}
Using (\ref{eq96}) we have:
\begin{eqnarray}
{\cal L}_1^{(n)}(z) \, L_2^+(w) &=& L_1^+(q^{c/2}z)^{t_1} \, 
(a_1^{-n})^{t_1} \, \left(\widehat 
R_{21}(q^{\frac{c}{2}}w/z)^{t_1}\right)^{-1} \, L_2^+(w) \, \widetilde 
L_1^-(z) \, \widehat R_{21}^*(q^{-\frac {c}{2}}w/z)^{t_1} \nonumber \\
&=& L_1^+(q^{c/2}z)^{t_1} \, (a_1^{-n})^{t_1} \, \left(\widehat 
R_{21}(q^{\frac{c}{2}}w/z)^{t_1}\right)^{-1} \, (a_1^{n})^{t_1} \, 
L_2^+(w) \, (a_1^{-n})^{t_1} \, \widetilde L_1^-(z) \, \widehat 
R_{21}^*(q^{-\frac {c}{2}}w/z)^{t_1} \,.  \nonumber \\
\label{eq98}
\end{eqnarray}
On the other hand, using the crossing-symmetry property we have:
\begin{eqnarray}
(a_1^{-n})^{t_1} \left(\widehat R_{21}(q^{\frac{c}{2}}w/z)^{t_1} 
\right)^{-1} (a_1^{n})^{t_1} &=& (a_1^{-n})^{t_1} \left(\widehat 
R_{21}(q^{\frac{c}{2}+N}w/z)^{-1} \right)^{t_1} (a_1^{n})^{t_1} 
\nonumber \\
&=& \left( a_1^{n}\widehat R_{21}(q^{\frac{c}{2}+N}w/z)^{-1} a_1^{-n} 
\right)^{t_1} \,.
\label{eq99}
\end{eqnarray}
We now apply $n$ times the following relation coming from
unitarity and quasi-periodicity:
\begin{equation}
\widehat R_{21} \left( z^{-1}(-p^{\frac{1}{2}}) \right)^{-1}=
\tau_N (q^{\frac{1}{2}}z) \tau_N (q^{\frac{1}{2}}z^{-1})
a_1 \widehat R_{21}(z^{-1})^{-1} a_1^{-1} \,.
\label {eq910}
\end{equation}

We see here the role of the quasi-periodicity operator in implementing
the general power of $-p^{\frac{1}{2}}$ in the $R$ matrix, leading to:
\begin{equation}
\widehat R_{21} \left(z^{-1}(-p^{\frac{1}{2}})^{n}\right)^{-1} = 
G_{N}(n,z) \, a_1^{n} \widehat R_{21}(z^{-1})^{-1} a_1^{-n} \,, \label{eq911}
\end{equation}
where:
\begin{eqnarray}
G_{N}(n,z) &=& \prod_{k=0}^{n-1} \tau_N \left[zq^{\frac{1}{2}} 
(-p^{\frac{1}{2}})^{-k} \right] \tau_N \left[z^{-1}q^{\frac{1}{2}} 
(-p^{\frac{1}{2}})^{k} \right] \quad \mbox{for $n>0$} \,, 
\nonumber \\
&& \label{eq911b} \\
G_{N}(n,z) &=& \prod_{k=1}^{|n|} \tau_N^{-1} \left[zq^{\frac{1}{2}} 
(-p^{\frac{1}{2}})^k \right] \tau_N^{-1} \left[z^{-1}q^{\frac{1}{2}} 
(-p^{\frac{1}{2}})^{-k}\right] \quad \mbox{for $n<0$} \,. 
\nonumber 
\end{eqnarray}
Applying (\ref{eq911}) to (\ref{eq99}) and using the equation 
(\ref{eq93}), we have:
\begin{equation}
(a_1^{-n})^{t_1} \left( \widehat R_{21}(q^{\frac{c}{2}}w/z)^{t_1} 
\right)^{-1} (a_1^{n})^{t_1} = G_{N}^{-1} 
\left(n,q^{\frac{c}{2}}(-p^{\frac{1}{2}})^{n} z/w \right) \left( 
\widehat R_{21}(q^{-\frac{c}{2}}w/z)^{-1} \right)^{t_1} \,.  
\label{eq912}
\end{equation}
Remark that from the definition of $\tau_N(x)$ (\ref{eq215}) and
the relation (\ref{eq93}), it follows that:
\begin{equation}
G_{N}^{-1}\left(n,q^{\frac{c}{2}}(-p^{\frac{1}{2}})^nx^{-1}\right) = 
F_{N}\left(n,q^{\frac{c}{2}}x\right) \,, \label{eq912b}
\end{equation}
where $F_{N}$ is given by (\ref{eq32}).  \\
Inserting (\ref{eq912}),(\ref{eq912b}) into (\ref{eq98}), we have:
\begin{equation}
{\cal L}_1^{(n)}(z) \, L_2^+(w) = 
F_{N}\left(n,q^{\frac{c}{2}}\frac{w}{z}\right) \, L_1^+(q^{c/2}z)^{t_1} 
\left( \widehat R_{21}(q^{-\frac{c}{2}}w/z)^{-1} \right)^{t_1} L_2^+(w) 
(a_1^{-n})^{t_1} \widetilde L_1^-(z) \widehat 
R^{*}_{21}(q^{-\frac{c}{2}}w/z)^{t_1} \,.
\label{eq913}
\end{equation}
Now we use equation (\ref{eq97}) to obtain:
\begin{equation}
{\cal L}_1^{(n)}(z) \, L_2^+(w) = 
F_{N}\left(n,q^{\frac{c}{2}}\frac{w}{z}\right) \, L^+_2(w) \left( 
\widehat R^{*}_{21}(q^{-\frac{c}{2}}w/z)^{-1} \right)^{t_1} 
L_1^+(q^{c/2}z)^{t_1} (a_1^{-n})^{t_1} \widetilde L_1^-(z) \widehat 
R^{*}_{21}(q^{-\frac{c}{2}}w/z)^{t_1} \,.  \label{eq914}
\end{equation}
\finproof

We are now able to state the following theorem:
\begin{thm}\label{thmsurf}
On the surfaces $\Sigma_{N,n}$, the operators $t(z)$ defined by 
(\ref{eq91}) satisfy the following exchange relations with the 
generators $L(w)$ of $\elp$:
\begin{equation}
t(z) \, L(w) = F_{N}\left(n,\frac{w}{z}\right) \, L(w) \, t(z) \,.
\label{eq94}
\end{equation}
\end{thm}

\medskip

\noindent \textbf{Proof:} We formulate the proof in terms of $L^+(w)$.  
One has
\begin{equation}
t(z) \, L_2^+(w) = {\mbox{Tr}}_{1} \Big[ L_1^+(zq^{c/2})^{t_1} \, 
(a_1^{-n})^{t_1} \, \widetilde L_1^-(z) \, L_2^+(w) \Big] \,.
\label{eq95bis}
\end{equation}
{From} eq. (\ref{eq914}), one obtains immediately:
\begin{equation}
t(z) \, L_2^+(w) = F_{N}\left(n,q^{\frac{c}{2}}\frac{w}{z}\right) \, 
L^+_2(w) {\mbox{Tr}}_{1} \left[ \left( \widehat 
R^{*}_{21}(q^{-\frac{c}{2}}w/z)^{-1} \right)^{t_1} L_1^+(q^{c/2}z)^{t_1} 
(a_1^{-n})^{t_1} \widetilde L_1^-(z) \widehat 
R^{*}_{21}(q^{-\frac{c}{2}}w/z)^{t_1} \right] \,.  \label{eq914bis}
\end{equation}
Using the very useful property:
\begin{equation}
{\mbox{Tr}}_{1} \left( R_{21} Q_1 R^\prime_{21} \right) = 
{\mbox{Tr}}_{1}\left( Q_1 {R^\prime_{21}}^{t_2} {R_{21}}^{t_2} 
\right)^{t_2} \,, \label{eq915}
\end{equation}
we get:
\begin{equation}
t(z) \, L_2^+(w) = F_{N}\left(n,q^{\frac{c}{2}}\frac{w}{z}\right) L_2^+(w) 
{\mbox{Tr}}_{1} \left[ L_1^+(q^{c/2}z)^{t_1} (a_1^{-n})^{t_1} \widetilde 
L_1^-(z) \widehat R^{*}_{21}(q^{-\frac{c}{2}}w/z)^{t} \left( \widehat 
R^{*}_{21}(q^{-\frac{c}{2}}w/z)^{-1} \right)^t \right]^{t_2} \,.  
\label{eq916}
\end{equation}
Now the product of the two $R$-matrices in 
(\ref{eq916}) vanishes altogether, leaving:
\begin{equation}
t(z) \, L^+(w) = F_{N}\left(n,q^\frac{c}{2}\frac{w}{z}\right) L^+(w) \, 
t(z) \,.  \label{eq917}
\end{equation}
{From} definitions (\ref{eq220}) Theorem \ref{thmsurf} follows then 
immediately.  
\finproof 

In particular since $a^N$ is proportionnal to the identity one recovers the 
exchange algebras in \cite{AFRS3} when $n = NM, \,\, M \in \ZZ$.

\medskip

A simple corollary of Theorem \ref{thmsurf} is the exchange relation 
between $t(z)$ and $t(w)$.  Indeed using the relations:
\begin{equation}
\widetilde L^-(z) = \left(L^-(z)^t\right)^{-1} \quad , \quad L^-(w) = 
aL^+(-p^{\frac{1}{2}}q^{-\frac{c}{2}}w)a^{-1} \,, \label{eq918}
\end{equation}
we derive from ({\ref{eq94}) the following result:

\begin{coro} 
On the surface $\Sigma_{N,n}$ the operators $t(z)$ defined by (\ref{eq91}) 
satisfy the following algebra:
\begin{equation}
t(z) \, t(w) = {\cal Y}_{N}(n,w/z) \, t(w)\, t(z) \,, \label{eq919}
\end{equation}
where ${\cal Y}_{N}$ is the function defined by (\ref{eq34}).
\end{coro}

\medskip

\noindent
\textbf{Remark:} For $n=-1$ equation (\ref{eq919}) reads as:
\begin{equation}
t(z) \, t(w) = t(w) \, t(z) \,, \label{eq921}
\end{equation}
that is, on the surface $q^{-c-N}=(-p^{\frac{1}{2}})^{-1}$, 
the operators $t(z)$ commute.  However they do not belong to the center 
of $\elp$, because the exchange factor of (\ref{eq94}) is different from 
unity.

\medskip

As when $c = -N$ this is an occurence of a one-step mechanism where one 
obtains directly a commuting algebra of operators with one single 
constraint on $p, q, c$.  However, contrary to that previous case, where 
the elliptic algebra and the quantum group shared this feature, this one 
is characteristic of the elliptic algebra structure, involving as it 
does the elliptic nome $p$.

\medskip

Exchange relations (\ref{eq919}) are to be understood as realizations of 
$Vir_{q,p}(sl(N))$ algebras in the framework of $\elp$.  This conclusion 
derives from the following results.

\begin{thm}\label{thmabel}
On the surface $\Sigma_{N,n}$, when $p = q^{Nh}$ with $h \in \ZZ 
\backslash \{ 0 \}$, the function ${\cal Y}(n,x)$ is equal to 1.  Hence 
$t(z)$ realizes an Abelian subalgebra in $\elp$.
\end{thm}

\medskip

\noindent 
\textbf{Proof:} The theorem \ref{thmabel} is easily proved using the 
periodicity properties of the $\Theta_{q^{2N}}$ functions.  
\finproof

\medskip

\begin{thm}\label{thmpois}
Setting $q^{Nh}=p^{1-\beta}$ for any integer $h \not = 0$, the 
$h$-labeled Poisson structure defined by:
\begin{equation}
\{ t(z) , t(w) \}^{(h)} = \lim_{\beta \rightarrow 0} \frac{1}{\beta}
\, \Big(t(z)t(w) - t(w)t(z)\Big)
\label{eq101}
\end{equation}
has the following expression:
\begin{equation}
\{ t(z) , t(w) \}^{(h)} = f_h(w/z) \, t(z) \, t(w)
\label{eq102}
\end{equation}
where
\begin{eqnarray}
f_h(x) &=&
2Nh \ln q \left[ \sum_{\ell \ge 0}
E(\frac{n}{2})(E(\frac{n}{2})+1)
\left( \frac{2x^2q^{2N\ell}}{1-x^2q^{2N\ell}}
- \frac{x^2q^{2N\ell+2}}{1-x^2q^{2N\ell+2}}
- \frac{x^2q^{2N\ell-2}}{1-x^2q^{2N\ell-2}} \right) \right.
\nonumber \\
&& + E(\frac{n+1}{2})^2 \left(
\frac{2x^2q^{2N\ell+N}}{1-x^2q^{2N\ell+N}}
- \frac{x^2q^{2N\ell+N+2}}{1-x^2q^{2N\ell+N+2}}
- \frac{x^2q^{2N\ell+N-2}}{1-x^2q^{2N\ell+N-2}} \right)
\nonumber \\
&& \left. - \frac{1}{2} E(\frac{n}{2})(E(\frac{n}{2})+1)
\left( \frac{2x^2}{1-x^2} - \frac{x^2q^2}{1-x^2q^2}
- \frac{x^2q^{-2}}{1-x^2q^{-2}} \right)
- (x \leftrightarrow x^{-1}) \right]
\qquad \mbox{for $h$ odd} \,, \nonumber\\
&& \label{eq103} \\
&=& Nh n(n+1) \ln q \left[ \sum_{\ell \ge 0} \left(
\frac{2x^2q^{2N\ell}}{1-x^2q^{2N\ell}}
- \frac{x^2q^{2N\ell+2}}{1-x^2q^{2N\ell+2}}
- \frac{x^2q^{2N\ell-2}}{1-x^2q^{2N\ell-2}} \right) \right.
\nonumber \\
&& \left. - \frac{1}{2} \left( \frac{2x^2}{1-x^2}
- \frac{x^2q^2}{1-x^2q^2} - \frac{x^2q^{-2}}{1-x^2q^{-2}} \right)
- (x \leftrightarrow x^{-1}) \right]
\qquad \mbox{for $h$ even} \,. \nonumber
\end{eqnarray}
Here the notation $E(m)$ means the integer part of the number $m$.
\end{thm}

\medskip

\noindent 
\textbf{Proof:} By direct calculation.
\finproof

\medskip

\noindent 
Formula (\ref{eq103}) is a trivial generalization of formulas (5.3) of 
\cite{AFRS3} provided that the formal substitution $NM \rightarrow n$ be 
done.  Therefore the discussion between Theorem 7 and Proposition 5 of 
\cite{AFRS3} is still valid and the Poisson bracket in the sector $k=0$ 
corresponding to (\ref{eq103}) is given, modulo the substitution 
indicated, by the formula in Proposition 5 of \cite{AFRS3}, which for 
even $h$ yields the algebra $Vir_{q}(sl(N))$ \cite{FR}.  Therefore we 
may conclude that (\ref{eq919}) realizes the exchange relations of 
\emph{bona fide} quantum $Vir_{q,p}(sl(N))$ algebras.

When $N=2$ the construction does not lead to new $q$-deformed Virasoro 
algebras.  However when $N \ge 3$, one obtain new exchange algebraic 
structures corresponding to the surfaces $\Sigma_{N,n}$ when $n \ne NM, 
\,\, M\in \ZZ$, and these structures will now be extended to complete 
$q$-deformed ${\cal W}_N$ algebras.

%%%%%%%%%%%%%%%%%%%%%%%%%%%%%%%%%%%%%%%%%%%%%%%%%%%%%%%%%%%%%%%%%%%%%%%
\section{Universal construction of ${\cal W}_{q,p}$ algebras}
\setcounter{equation}{0}

Extending the construction to higher spin currents for $N \ge 3$ compels 
to use the direct exchange relation in Lemma 1 instead of the quadratic 
intertwined exchange relation in Theorem 2, for which no generalizations 
to higher powers of ${\mathfrak T}$ exist.

\subsection{Quantum ${\cal W}_{q,p}$ algebras}

\begin{thm}\label{thmwl}
We define the operators $w_{s}(z)$ ($s = 1, \dots, N-1$) by:
\begin{equation}
w_{s}(z) \equiv \mbox{Tr} \left[ \left( \prod_{1 \le i \le 
s}^{{\textstyle \curvearrowleft}} \prod_{j>i} P_{ij} \right) \prod_{1 
\le i \le s}^{{\textstyle \curvearrowleft}} \left( 
{\cal L}_{i}^{(n)}(z_i) \prod_{j>i} \widehat 
R_{ij}^{*}(q^{-N}z_i/z_j)^{t_it_j} \right) \right] \,.
\label{eq41}
\end{equation}
where 
\begin{equation}
{\cal L}_{i}^{(n)}(z) \equiv (a^{-n} \, 
L^+_{i}(q^{\frac{c}{2}}z))^{t_i} \, \widetilde L^-_{i}(z) \equiv 
\underbrace{\II \otimes \dots \otimes \II}_{i-1} \otimes (a^{-n} \, 
L^+(q^{\frac{c}{2}}z))^{t} \, \widetilde L^-(z) \otimes \underbrace{\II 
\otimes \dots \II}_{s-i}
\end{equation}
with $n \in \ZZ$, $z_i = zq^{i-\frac{s+1}{2}}$, and $P_{ij}$ is the 
permutation operator between the spaces $i$ and $j$ including
the spectral parameters.  \\
On the surface $\Sigma_{N,n}$ defined by $(-p^{\frac{1}{2}})^{n} = 
q^{-c-N}$, the operators $w_{s}(z)$ realize an exchange algebra with 
the generators $L(w)$ of $\elp$:
\begin{equation}
w_{s}(z) \, L(w) = F_{N}^{(s)}\Big(n,\frac{w}{z}\Big) \, L(w) \, 
w_{s}(z) \,,
\label{eq42}
\end{equation}
where
\begin{equation}
F_{N}^{(s)}\Big(n,\frac{w}{z}\Big) = \prod_{i=1}^{s} F_{N}\Big( 
n,\frac{w}{z_i} \Big) \,.
\end{equation}
\end{thm}

\medskip

\noindent \textbf{Proof:} For simplicity, we will only prove the theorem 
for $w_{2}(z)$ and $w_{3}(z)$ (the proof for $w_{1}(z) \equiv t(z)$ has 
been done in \cite{AFRS2,AFRS3}, see theorem \ref{thmvir} above).  \\
The proof is based on the exchange relation (\ref{eq914}) between ${\cal 
L}^{(n)}(z)$ and $L^+(w)$ on the surface $\Sigma_{N,n}$ defined by 
$(-p^{\frac{1}{2}})^{n} = q^{-c-N}$:
\begin{equation}
{\cal L}_i^{(n)}(z) \, L^+_\alpha(w) = 
F_{N}\left(n,q^{\frac{c}{2}}\frac{w}{z}\right) \, L_\alpha^+(w) \, 
\left( \widehat R^{*}_{\alpha i}(q^{-\frac{c}{2}}w/z)^{-1} \right)^{t_i} 
\, {\cal L}_i^{(n)}(z) \, \widehat R^{*}_{\alpha 
i}(q^{-\frac{c}{2}}w/z)^{t_i} \,.  \label{eq410}
\end{equation}

\bigskip

\noindent Consider the operator $w_{2}(z)$.  By definition (with $z_1 = 
zq^{-\half}$ and $z_2 = zq^{\half}$):
\begin{equation}
w_{2}(z) = \mbox{Tr} \Big[ P_{12} \, {\cal L}_{1}^{(n)}(z_1) \, \widehat 
R_{12}^{*}(q^{-N}z_1/z_2)^{t_1t_2} \, {\cal L}_{2}^{(n)}(z_2) \Big] \,.
\end{equation}
{From} the exchange relation (\ref{eq410}) between ${\cal L}^{(n)}(z)$ 
and $L^+(w)$, one gets immediately:
\begin{eqnarray}
w_{2}(z) \, L^+_{\alpha}(w) &=& \mbox{Tr}_{12} \Big[ P_{12} \, {\cal 
L}_{1}^{(n)}(z_1) \, \widehat R_{12}^{*}(q^{-N}z_1/z_2)^{t_1t_2} \, {\cal 
L}_{2}^{(n)}(z_2) \Big] \, L^+_{\alpha}(w) \nonumber \\
&=& F_{N}^{(2)}\Big(n,q^{\frac{c}{2}}\frac{w}{z}\Big) \, L^+_{\alpha}(w) 
\, \mbox{Tr}_{12} \Big[ P_{12} \left( \widehat R_{\alpha 
1}^*(q^{-\frac{c}{2}}w/z_1)^{-1} \right)^{t_1} \, {\cal L}_{1}^{(n)}(z_1) \, 
{\widehat R_{\alpha 1}^*(q^{-\frac{c}{2}}w/z_1)}^{t_1}\nonumber \\
&& \widehat R_{12}^*(q^{-N}z_1/z_2)^{t_1t_2} \, \left( 
\widehat R_{\alpha 2}^*(q^{-\frac{c}{2}}w/z_2)^{-1} \right)^{t_2} 
\, {\cal L}_{2}^{(n)}(z_2) \, 
{\widehat R_{\alpha 2}^*(q^{-\frac{c}{2}}w/z_2)}^{t_2} \Big] \,,
\label{eq412}
\end{eqnarray}
where $\displaystyle F_{N}^{(2)}\Big(n,q^{\frac{c}{2}}\frac{w}{z}\Big)$ 
is given by:
\begin{equation}
F_{N}^{(2)}\Big(n,q^{\frac{c}{2}}\frac{w}{z}\Big) = F_{N}\Big( 
n,q^{\frac{c}{2}}\frac{w}{z_1} \Big) \, F_{N}\Big( 
n,q^{\frac{c}{2}}\frac{w}{z_2} \Big) \,.
\end{equation}
In order to reorganize the $R$ matrices in (\ref{eq412}), one uses the 
Yang--Baxter equation for the matrix $\widehat R^*$ (equation 
(\ref{eq413}) is a consequence of (\ref{eq29}), the normalization factor 
entering in the definition (\ref{eq214}) of the $\widehat R$ matrices 
being the same in the l.h.s.  and in the r.h.s.  of (\ref{eq413})):
\begin{equation}
\widehat R_{\alpha 1}^*(x_1) \, \widehat R_{\alpha 2}^*(x_2) \, 
\widehat R_{12}^*(x_2/x_1) = \widehat R_{12}^*(x_2/x_1) \, 
\widehat R_{\alpha 2}^*(x_2) \, \widehat R_{\alpha 1}^*(x_1) \,,
\label{eq413}
\end{equation}
{from} which it follows (with a shift $x_2 \rightarrow q^{-N} x_2$)
\begin{equation}
\widehat R_{\alpha 1}^*(x_1)^{t_1} \, 
\widehat R_{12}^*(q^{-N}x_2/x_1)^{t_1t_2} \, 
\left( \widehat R_{\alpha 2}^*(x_2)^{-1} \right)^{t_2} = 
\left( \widehat R_{\alpha 2}^*(x_2)^{-1} \right)^{t_2} \, 
\widehat R_{12}^*(q^{-N}x_2/x_1)^{t_1t_2} \, 
\widehat R_{\alpha 1}^*(x_1)^{t_1} \,.
\label{eq414}
\end{equation}
Therefore, one obtains
\begin{eqnarray}
w_{2}(z) \, L^+_{\alpha}(w) &=& 
F_{N}^{(2)}\Big(n,q^{\frac{c}{2}}\frac{w}{z}\Big) \, 
L^+_{\alpha}(w) \, \mbox{Tr}_{12} \Big[ P_{12} \, \left( \widehat 
R_{\alpha 1}^*(q^{-\frac{c}{2}}w/z_1)^{-1} \right)^{t_1} \, {\cal 
L}_{1}^{(n)}(z_1) \, \left( \widehat R_{\alpha 
2}^*(q^{-\frac{c}{2}}w/z_2)^{-1} \right)^{t_2} \nonumber \\
&& \widehat R_{12}^*(q^{-N}z_1/z_2)^{t_1t_2} \, 
{\widehat R_{\alpha 1}^*(q^{-\frac{c}{2}}w/z_1)}^{t_1} 
\, {\cal L}_{2}^{(n)}(z_2) \, {\widehat 
R_{\alpha 2}^*(q^{-\frac{c}{2}}w/z_2)}^{t_2} \Big] \nonumber \\
&=& F_{N}^{(2)}\Big(n,q^{\frac{c}{2}}\frac{w}{z}\Big) \, L^+_{\alpha}(w) \, 
\mbox{Tr}_{12} \Big[ P_{12} \, \left( \widehat R_{\alpha 
1}^*(q^{-\frac{c}{2}}w/z_1)^{-1} \right)^{t_1} \, \left( \widehat 
R_{\alpha 2}^*(q^{-\frac{c}{2}}w/z_2)^{-1} \right)^{t_2} \nonumber \\
&& {\cal L}_{1}^{(n)}(z_1) \, \widehat R_{12}^*(q^{-N}z_1/z_2)^{t_1t_2} 
\, {\cal L}_{2}^{(n)}(z_2) \, {\widehat R_{\alpha 
1}^*(q^{-\frac{c}{2}}w/z_1)}^{t_1} \, {\widehat R_{\alpha 
2}^*(q^{-\frac{c}{2}}w/z_2)}^{t_2} \Big] \nonumber \\ 
&=& F_{N}^{(2)}\Big(n,q^{\frac{c}{2}}\frac{w}{z}\Big) \, L^+_{\alpha}(w) \, 
\mbox{Tr}_{12} \Big[ \left( \widehat R_{\alpha 
2}^*(q^{-\frac{c}{2}}w/z_2)^{-1} \right)^{t_2} \, \left( \widehat 
R_{\alpha 1}^*(q^{-\frac{c}{2}}w/z_1)^{-1} \right)^{t_1} \, P_{12} 
\nonumber \\
&& {\cal L}_{1}^{(n)}(z_1) \, \widehat R_{12}^*(q^{-N}z_1/z_2)^{t_1t_2} 
\, {\cal L}_{2}^{(n)}(z_2) \, {\widehat R_{\alpha 
1}^*(q^{-\frac{c}{2}}w/z_1)}^{t_1} \, {\widehat R_{\alpha 
2}^*(q^{-\frac{c}{2}}w/z_2)}^{t_2} \Big] \,. 
\label{eq415}
\end{eqnarray}
the last equality being obtained by action of the permutation operator 
$P_{12}$.  \\

One then uses the fact that under a trace over the space $\beta$, for 
any $c$-number matrices $A_{\alpha \beta}$ and $B_{\alpha \beta}$ and 
operator matrix $Q_{\alpha}$, one has
\begin{equation}
\mbox{Tr}_{\beta} \Big[ A_{\alpha \beta} Q_{\beta} B_{\alpha \beta} 
\Big] = \mbox{Tr}_{\beta} \Big[ Q_{\beta} \Big((B_{\alpha 
\beta})^{t_\alpha} (A_{\alpha \beta})^{t_\alpha}\Big)^{t_\alpha} \Big] 
\,,
\label{eq416}
\end{equation}
Applying (\ref{eq416}) to $\beta \equiv 1 \otimes 2, A_{\alpha \beta} 
\equiv R_{\alpha 2} R_{\alpha 1}$ and $B_{\alpha \beta} \equiv 
R'_{\alpha 1} R'_{\alpha 2}$, one gets
\begin{eqnarray}
w_{2}(z) \, L^+_{\alpha}(w) &=& 
F_{N}^{(2)}\Big(n,q^{\frac{c}{2}}\frac{w}{z}\Big) \, 
L^+_{\alpha}(w) \, \mbox{Tr}_{12} \Big[ P_{12} \, {\cal L}_{1}^{(n)}(z_1) \, 
\widehat R_{12}^*(q^{-N}z_1/z_2)^{t_1t_2} \, {\cal L}_{2}^{(n)}(z_2) \, 
\Big({\widehat R_{\alpha 2}^*(q^{-\frac{c}{2}}w/z_2)}^{t_2t_\alpha}
\nonumber \\
&& {\widehat R_{\alpha 1}^*(q^{-\frac{c}{2}}w/z_1)}^{t_1t_\alpha} \, 
\left( \widehat R_{\alpha 1}^*(q^{-\frac{c}{2}}w/z_1)^{-1} 
\right)^{t_1t_\alpha} \, \left( \widehat R_{\alpha 
2}^*(q^{-\frac{c}{2}}w/z_2)^{-1} \right)^{t_2t_\alpha}\Big)^{t_\alpha} 
\Big] \nonumber \\
&=& F_{N}^{(2)}\Big(n,q^{\frac{c}{2}}\frac{w}{z}\Big) \, L^+_{\alpha}(w) \, 
\mbox{Tr}_{12} \Big[ P_{12} \, {\cal L}_{1}^{(n)}(z_1) \, \widehat 
R_{12}^*(q^{-N}z_1/z_2)^{t_1t_2} \, {\cal L}_{2}^{(n)}(z_2) \, 
\Big({\widehat R_{\alpha 2}^*(q^{-\frac{c}{2}}w/z_2)}^{t_2t_\alpha} 
\nonumber \\
&& {\widehat R_{\alpha 1}^*(q^{-\frac{c}{2}}w/z_1)}^{t_1t_\alpha} \, 
\left( \widehat R_{\alpha 1}^*(q^{-\frac{c}{2}}w/z_1)^{t_1t_\alpha} 
\right)^{-1} \, \left( \widehat R_{\alpha 
2}^*(q^{-\frac{c}{2}}w/z_2)^{t_2t_\alpha} \right)^{-1}\Big)^{t_\alpha} 
\Big] \nonumber \\
&=& F_{N}^{(2)}\Big(n,q^{\frac{c}{2}}\frac{w}{z}\Big) \, L^+_{\alpha}(w) \, 
\mbox{Tr}_{12} \Big[ P_{12} \, {\cal L}_{1}^{(n)}(z_1) \, \widehat 
R_{12}^*(q^{-N}z_1/z_2)^{t_1t_2} \, {\cal L}_{2}^{(n)}(z_2) \Big] \,.
\end{eqnarray}
It follows that
\begin{equation}
w_{2}(z) \, L^+_{\alpha}(w) = 
F_{N}^{(2)}\Big(n,q^{\frac{c}{2}}\frac{w}{z}\Big) \, L^+_{\alpha}(w) \, 
w_{2}(z) \,.
\end{equation}
Recalling the fact that $L^+_{\alpha}(w) = 
L_{\alpha}(q^{\frac{c}{2}}w)$, one gets the desired result.

\bigskip

\noindent Consider now the case of $w_{3}(z)$.  By definition (with $z_1 
= zq^{-1}$, $z_2 = z$ and $z_3 = zq$):
\begin{equation}
w_{3}(z) = \mbox{Tr} \Big[ P_{12} \, P_{13} \, P_{23} \, {\cal 
L}_{1}^{(n)}(z_1) \, \widehat R_{12}^{*}(q^{-N}z_1/z_2)^{t_1t_2} \, \widehat 
R_{13}^{*}(q^{-N}z_1/z_3)^{t_1t_3} \, {\cal L}_{2}^{(n)}(z_2) \, \widehat 
R_{23}^{*}(q^{-N}z_2/z_3)^{t_2t_3} \, {\cal L}_{3}^{(n)}(z_3) \Big] \,.
\end{equation}
{From} the exchange relation (\ref{eq410}) between ${\cal L}^{(n)}(z)$ and 
$L^+(w)$, one gets:
\begin{eqnarray}
w_{3}(z) \, L^+_{\alpha}(w) &=& \mbox{Tr}_{123} \Big[ P_{12} \, P_{13} 
\, P_{23} \, {\cal L}_{1}^{(n)}(z_1) \, \widehat 
R_{12}^{*}(q^{-N}z_1/z_2)^{t_1t_2} \, \widehat 
R_{13}^{*}(q^{-N}z_1/z_3)^{t_1t_3} \nonumber \\
&& {\cal L}_{2}^{(n)}(z_2) \, \widehat R_{23}^{*}(q^{-N}z_2/z_3)^{t_2t_3} \, 
{\cal L}_{3}^{(n)}(z_3) \Big] \, L^+_{\alpha}(w) \nonumber \\
&=& F_{N}^{(3)}\Big(n,q^{\frac{c}{2}}\frac{w}{z}\Big) \, L^+_{\alpha}(w) \, 
\mbox{Tr}_{123} \Big[ P_{12} \, P_{13} \, P_{23} \, \left( \widehat 
R_{\alpha 1}^*(q^{-\frac{c}{2}}w/z_1)^{-1} \right)^{t_1} \, {\cal 
L}_{1}^{(n)}(z_1) \, {\widehat R_{\alpha 1}^*(q^{-\frac{c}{2}}w/z_1)}^{t_1} 
\nonumber \\
&& \widehat R_{12}^{*}(q^{-N}z_1/z_2)^{t_1t_2} \, \widehat 
R_{13}^{*}(q^{-N}z_1/z_3)^{t_1t_3} \, \left( \widehat R_{\alpha 
2}^*(q^{-\frac{c}{2}}w/z_2)^{-1} \right)^{t_2} \, {\cal 
L}_{2}^{(n)}(z_2) \, {\widehat R_{\alpha 
2}^*(q^{-\frac{c}{2}}w/z_2)}^{t_2} \nonumber \\
&& \widehat R_{23}^{*}(q^{-N}z_2/z_3)^{t_2t_3} \, \left( \widehat 
R_{\alpha 3}^*(q^{-\frac{c}{2}}w/z_3)^{-1} \right)^{t_3} \, {\cal 
L}_{3}^{(n)}(z_3) \,{\widehat R_{\alpha 
3}^*(q^{-\frac{c}{2}}w/z_3)}^{t_3} \Big] \,,
\end{eqnarray}
where $\displaystyle F_{N}^{(3)}\Big(n,q^{\frac{c}{2}}\frac{w}{z}\Big)$ 
is given by:
\begin{equation}
F_{N}^{(3)}\Big(n,q^{\frac{c}{2}}\frac{w}{z}\Big) = F_{N}\Big( 
n,q^{\frac{c}{2}}\frac{w}{z_1} \Big) \, F_{N}\Big( 
n,q^{\frac{c}{2}}\frac{w}{z_2} \Big) \, F_{N}\Big( 
n,q^{\frac{c}{2}}\frac{w}{z_3} \Big) \,.
\end{equation}
Applying three times the Yang--Baxter equation (\ref{eq414}), one obtains:
\begin{eqnarray}
w_{3}(z) \, L^+_{\alpha}(w) &=& 
F_{N}^{(3)}\Big(n,q^{\frac{c}{2}}\frac{w}{z}\Big) \, L^+_{\alpha}(w) \, 
\mbox{Tr}_{123} \Big[ P_{12} \, P_{13} \, P_{23} \, \left( \widehat 
R_{\alpha 1}^*(q^{-\frac{c}{2}}w/z_1)^{-1} \right)^{t_1} \, \left( 
\widehat R_{\alpha 2}^*(q^{-\frac{c}{2}}w/z_2)^{-1} \right)^{t_2} 
\nonumber \\
&& \left( \widehat R_{\alpha 3}^*(q^{-\frac{c}{2}}w/z_3)^{-1} 
\right)^{t_3} \, {\cal L}_{1}^{(n)}(z_1) \, \widehat 
R_{12}^{*}(q^{-N}z_1/z_2)^{t_1t_2} \, \widehat 
R_{13}^{*}(q^{-N}z_1/z_3)^{t_1t_3} \, {\cal L}_{2}^{(n)}(z_2) \nonumber 
\\
&& \widehat R_{23}^{*}(q^{-N}z_2/z_3)^{t_2t_3} \, {\cal L}_{3}^{(n)}(z_3) \, 
{\widehat R_{\alpha 1}^*(q^{-\frac{c}{2}}w/z_1)}^{t_1} \, {\widehat 
R_{\alpha 2}^*(q^{-\frac{c}{2}}w/z_2)}^{t_2} \, {\widehat R_{\alpha 
3}^*(q^{-\frac{c}{2}}w/z_3)}^{t_3} \Big] \,.  \nonumber \\
\label{eq422}
\end{eqnarray}
Finally, after action of the permutation operators and using (\ref{eq416}) 
applied to $\beta \equiv 1 \otimes 2 \otimes 3$, the $R$ matrices 
$\widehat R_{\alpha i}^*$ in (\ref{eq422}) simplify.  It follows that
\begin{equation}
w_{3}(z) \, L^+_{\alpha}(w) = 
F_{N}^{(3)}\Big(n,q^{\frac{c}{2}}\frac{w}{z}\Big) \, L^+_{\alpha}(w) \, 
w_{3}(z) \,.
\end{equation}
Again, one gets the desired result since $L^+_{\alpha}(w) = 
L_{\alpha}(q^{\frac{c}{2}}w)$.

\bigskip

\noindent The proof for a generic operator $w_{s}(z)$ is obtained by 
using the basic exchange relation (\ref{eq410}) between ${\cal 
L}^{(n)}(z)$ and $L^+(w)$ and applying $\sfrac{1}{2} s(s-1)$ times the 
Yang--Baxter equation in the form (\ref{eq414}).  Successive uses of 
this relation yields an expression involving: the group of permutation 
operators; the product of all $R$ matrices appearing at the left of the 
${\cal L}^{(n)}$ operator in Lemma 1; the terms of the monomial 
$w_{s}(z)$; finally the product of all $R$ matrices appearing at the 
right of the ${\cal L}^{(n)}$ operator in Lemma 1.  Commutation of the 
permutation operators with the ``left'' $R$ matrices then brings this 
group of $R$ matrices in position to use the transposition procedure 
(\ref{eq416}), and precisely rearranges the indices of these $R$ 
matrices in the exact way required to cancel the two ``left'' and 
``right'' groups of exchange-generated $R$ matrices once the 
generalization of the procedure (\ref{eq416}) is used.  Therefore the 
exchange relation becomes an exchange algebra with purely scalar 
functional structure coefficients.  
\finproof

\bigskip

An immediate consequence is therefore the exchange relation between the 
operators $w_{s}(z)$ and the generators ${\cal L}^{(n)}(w)$:
\begin{equation}
w_{s}(z) \, {\cal L}^{(n)}(w) = \frac{\displaystyle 
F_{N}^{(s)}\Big(n,q^c\frac{w}{z}\Big)} {\displaystyle 
F_{N}^{(s)}\Big(n,-p^{\frac{1}{2}}\frac{w}{z}\Big)} \, {\cal L}^{(n)}(w) 
\, w_{s}(z) = \prod_{i=1}^s {\cal Y}_{N}\Big(n,\frac{w}{z_{i}}\Big) \, 
{\cal L}^{(n)}(w) \, w_{s}(z) \,.  \label{eq425}
\end{equation}
where ${\cal Y}_{N}$ is the function defined by (\ref{eq34}) and 
$z_i = zq^{i-\frac{s+1}{2}}$.

\medskip

It is now immediate to derive the exchange algebra among the operators 
$w_{s}(z)$.  One gets the following theorem:
\begin{thm}\label{thmww}
On the surface $\Sigma_{N,n}$, the operators $w_{s}(z)$ realize an 
exchange algebra
\begin{equation}
w_i(z) \, w_j(w) = \prod_{u=-\frac{i-1}{2}}^{\frac{i-1}{2}}
\prod_{v=-\frac{j-1}{2}}^{\frac{j-1}{2}}{\cal Y}_{N} 
\left(n,q^{v-u}\frac{w}{z}\right) w_j(w) \, w_i(z) \,.
\end{equation}
\end{thm}
The proof is obvious and follows immediately from definition \ref{eq41} 
and equation (\ref{eq425}).
\finproof

\medskip

\noindent \textbf{Remarks}: \\
1) The critical level $c=-N$ can be seen as a limiting case of the 
relation $(-p^{\frac{1}{2}})^{n} = q^{-c-N}$ by taking $n=0$, $p$ and 
$q$ having arbitrary generic values ($\vert q \vert < 1$, $\vert p \vert 
< 1$).  In this limiting case, it is easy to note that the factor 
$G_{N}(0,x)$ is equal to 1 (see e.g.  (\ref{eq911})) and hence $F(0,x) = 
1$.  Therefore, at the critical level, the operators $w_{s}(z)$ provide 
\emph{a set of commuting quantities belonging to the center of the 
elliptic quantum $\elp$ algebra}.

\medskip

\noindent 
2) When one chooses $n = NM$ for $M \in \ZZ$, one recovers the structure 
functions of the $q$-${\cal W}_{N}$ algebras constructed in 
\cite{AFRS3}.  General values of $n$ lead to original structure 
functions.  This time however the construction does not go through the 
two steps of first constructing the trace $s(z)$ and then shift-multiply 
the same \emph{derived} abstract generator to get the full $q$-${\cal 
W}_{N}$ algebra: the trace procedure \emph{and} shifted multiplication 
are applied in one single stroke to the \emph{original} elliptic algebra 
generators, hence our denomination of this as a ``universal'' procedure.

\subsection{Classical limit}

\begin{thm}\label{thmclass}
On the surface $\Sigma_{N,n}$, when a additional relation $p = q^{Nh}$ 
with $h \in \ZZ \backslash \{0\}$ is imposed, the function ${\cal Y}$ is 
equal to 1.  Hence, the operators $w_{s}(z)$ realize an Abelian 
subalgebra in $\elp$.
\end{thm}

\medskip

\noindent \textbf{Proof}: The proof is straightforward by using the 
explicit expression of the function $F_{s}(M,x)$ and the periodicity 
properties of the $\Theta_{q^{2N}}$ functions.  
\finproof

\medskip

\noindent \textbf{Remark}:
Although the $w_{s}(z)$ realize an Abelian algebra, they do not belong 
to the center of $\elp$, in contrast to the case of the critical level.

The theorem \ref{thmclass} allows us to define Poisson structures on the 
corresponding Abelian subalgebras.  As usual, they are obtained as 
limits of the exchange algebras when $p = q^{Nh}$ with $h \in \ZZ 
\backslash \{0\}$.

\begin{thm}\label{thmclass2}
Setting $q^{Nh}=p^{1-\beta}$ for any integer $h \not = 0$, the
$h$-labeled Poisson structure defined by:
\begin{equation}
\{ w_{i}(z) , w_{j}(w) \}^{(h)} = \lim_{\beta \rightarrow 0} \frac{1}{\beta}
\, \Big( w_{i}(z) \, w_{j}(w) - w_{j}(w) \, w_{i}(z) \Big)
\label{eq430}
\end{equation}
has the following expression:
\begin{equation}
\{ w_{i}(z),w_{j}(w) \} = \sum_{u=-(i-1)/2}^{(i-1)/2}
\sum_{v=-(j-1)/2}^{(j-1)/2}
f_h \Big(q^{v-u}\frac{w}{z}\Big) \, w_{i}(z) \, w_{j}(w) \,,
\label{eq431}
\end{equation}
where
$f_h(x)$ is given by (\ref{eq103}).
\end{thm}

\medskip

One recovers here the same Poisson algebra as in \cite{AFRS3}, 
identifying therefore the exchange algebras in Theorem 7 as new ${\cal 
W}_{p,q}[sl(N)]$ algebras.  The main issue now is to obtain explicit 
realizations of these algebras.  Curiously enough the only known 
($q$-boson type) explicit realization achieved in \cite{FR} corresponds 
to a ${\cal W}_N$ algebra which does not belong to the set constructed 
here.  We hope to adress this issue in the next future.

\end{document}